\documentclass{amsart}
\usepackage{amssymb}
\usepackage{amsmath}
\usepackage{amsfonts}

\newtheorem{theorem}{Theorem}
\theoremstyle{plain}

\newtheorem{corollary}{Corollary}

\newtheorem{lemma}{Lemma}

\newtheorem{remark}{Remark}

\numberwithin{equation}{section}
\input{tcilatex}

\begin{document}
\title[Sampling Theory for Sturm-Liouville Problem ]{Sampling Theory for
Sturm-Liouville Problem with Boundary and Transmission Conditions Dependent
on Eigenparameter\ \ \ }
\author{Fatma H\i ra}
\address[F. H\i ra and N. Alt\i n\i \c{s}\i k]{Ondokuz May\i s University,
Arts and Science Faculty, 55139, Samsun, Turkey }
\email[Corresponding author: F. H\i ra]{fatma.hira@omu.edu.tr}
\author{Nihat Alt\i n\i \c{s}\i k}
\email[N. Alt\i n\i \c{s}\i k]{anihat@omu.edu.tr}

\begin{abstract}
In this paper, we investigate the sampling analysis associated with
discontinuous Sturm-Liouville problem which has transmission conditions at
two points of discontinuity also contains an eigenparameter in a boundary
condition and two transmission conditions. We establish briefly spectral
properties of the problem and then we prove the sampling theorem associated
with the problem.
\end{abstract}

\maketitle

\section{Introduction}

Let $\sigma >0,$ and denote by $B_{\sigma }^{2},$ the Paley-Wiener space of
all entire functions $f$ of exponential type with band width at most $\sigma 
$ which are $L^{2}\left( 
%TCIMACRO{\U{211d} }%
%BeginExpansion
\mathbb{R}
%EndExpansion
\right) $ functions when restricted to $%
%TCIMACRO{\U{211d} }%
%BeginExpansion
\mathbb{R}
%EndExpansion
.$ This space is characterized by the following relation which is given by
Paley and Wiener $\left[ 1\right] ,$ known as the Paley-Wiener theorem:%
\begin{equation}
f\left( t\right) \in B_{\sigma }^{2}\Leftrightarrow f\left( t\right) =\frac{1%
}{\sqrt{2\pi }}\int\limits_{-\sigma }^{\sigma }e^{ixt}g\left( x\right) dx,%
\text{ \ \ \ \ \ }g\in L^{2}\left( -\sigma ,\sigma \right) .  \tag{1.1}
\end{equation}%
\qquad \qquad \qquad

The Whittaker-Kotel'nikov-Shannon (WKS) sampling theorem states that if $%
f\left( t\right) \in B_{\sigma }^{2},$ then it is completely determined by
its values at the points $t_{k}=\frac{k\pi }{\sigma },k\in 
%TCIMACRO{\U{2124} }%
%BeginExpansion
\mathbb{Z}
%EndExpansion
$ and can be reconstructed by means of the formula%
\begin{equation}
f\left( t\right) =\sum\limits_{k=-\infty }^{\infty }f\left( t_{k}\right)
\sin c\left( \sigma t-k\pi \right) ,\text{ \ \ }t\in 
%TCIMACRO{\U{2102} }%
%BeginExpansion
\mathbb{C}
%EndExpansion
,  \tag{1.2}
\end{equation}%
where%
\begin{equation}
\sin c\left( \sigma t-k\pi \right) =\left\{ \QATOP{\frac{\sin \left( \sigma
t-k\pi \right) }{\sigma t-k\pi },\text{ \ \ \ \ }t\neq \frac{k\pi }{\sigma },%
}{1,\text{ \ \ \ \ \ \ \ \ \ \ \ \ \ \ \ }t=\frac{k\pi }{\sigma }.}\right. 
\tag{1.3}
\end{equation}%
The series (1.2) is absolutely and uniformly convergent on compact subsets
of $%
%TCIMACRO{\U{2102} }%
%BeginExpansion
\mathbb{C}
%EndExpansion
,$ uniformly convergent on $%
%TCIMACRO{\U{211d} }%
%BeginExpansion
\mathbb{R}
%EndExpansion
,($see $\left[ 2,3\right] ).$

One of the important generalizations of the WKS sampling theorem is the
Paley-Wiener-Levinson theorem which can be stated as follows: Let $\left\{
t_{k}\right\} _{k\in 
%TCIMACRO{\U{2124} }%
%BeginExpansion
\mathbb{Z}
%EndExpansion
}$ be a sequence of real numbers satisfying%
\begin{equation}
D:=\sup_{k\in 
%TCIMACRO{\U{2124} }%
%BeginExpansion
\mathbb{Z}
%EndExpansion
}\left\vert t_{k}-\frac{k\pi }{\sigma }\right\vert <\frac{\pi }{4\sigma }, 
\tag{1.4}
\end{equation}%
and let 
\begin{equation}
G\left( t\right) :=\left( t-t_{0}\right) \prod\limits_{k=1}^{\infty }\left(
1-\frac{t}{t_{k}}\right) \left( 1-\frac{t}{t_{-k}}\right) .  \tag{1.5}
\end{equation}%
Then for any $f\left( t\right) \in B_{\sigma }^{2},$%
\begin{equation}
f\left( t\right) =\sum\limits_{k=-\infty }^{\infty }f\left( t_{k}\right) 
\frac{G\left( t\right) }{G^{\prime }\left( t_{k}\right) \left(
t-t_{k}\right) },\text{ \ \ \ }t\in 
%TCIMACRO{\U{2102} }%
%BeginExpansion
\mathbb{C}
%EndExpansion
.  \tag{1.6}
\end{equation}%
The series (1.6) converges uniformly an compact subsets of $%
%TCIMACRO{\U{2102} }%
%BeginExpansion
\mathbb{C}
%EndExpansion
,($see $\left[ 2,4\right] ).$

Series of the form (1.6) is called a Lagrange-type interpolation.

Another important generalization of the WKS sampling theorem is the theorem
of Kramer $\left[ 5\right] ,$ can be stated as follows: Let $I$ be a finite
closed interval, $K\left( x,t\right) $ a function continuous in $t$~such
that $K\left( x,t\right) \in L^{2}\left( I\right) ~$for all $t\in 
%TCIMACRO{\U{2102} }%
%BeginExpansion
\mathbb{C}
%EndExpansion
,$ and let $\left\{ t_{k}\right\} _{k\in 
%TCIMACRO{\U{2124} }%
%BeginExpansion
\mathbb{Z}
%EndExpansion
}$ be a sequence of real numbers such that $\left\{ K\left( x,t_{k}\right)
\right\} _{k\in 
%TCIMACRO{\U{2124} }%
%BeginExpansion
\mathbb{Z}
%EndExpansion
}$ is a complete ortogonal set in $L^{2}\left( I\right) .$ Suppose that%
\begin{equation}
f\left( t\right) =\int\limits_{I}K\left( x,t\right) g\left( x\right) dx,%
\text{ \ \ \ \ }g\in L^{2}\left( I\right) .  \tag{1.7}
\end{equation}%
Then%
\begin{equation}
f\left( t\right) =\sum\limits_{k=-\infty }^{\infty }f\left( t_{k}\right) 
\frac{\int\limits_{I}K\left( x,t\right) \overline{K\left( x,t_{k}\right) }dx%
}{\left\Vert K\left( x,t_{k}\right) \right\Vert _{L^{2}\left( I\right) }^{2}}%
.  \tag{1.8}
\end{equation}

Generalization of WKS sampling theorem has been investigated extensively,
see $\left[ 6-9\right] .$ Sampling theorems associated with Sturm-Liouville
problems were investigated in $\left[ 10-15\right] $. Also $\left[ 16,17%
\right] $ and $\left[ 18,19\right] $ are the example works with an
eigenparameter in the boundary conditions in direction of sampling analysis
associated with continuous and discontinuous eigenproblems, respectively. In 
$\left[ 20\right] ,$ the authors discussed the situation of deriving a
sampling theorem of Kramer type when the kernels are discontinuous. They
introduced the discontinuous Sturm-Liouville problems studied by Kobayashi $%
\left[ 21\right] $ to Sturm-Liouville problems with eigenfunctions having
two symmetrically located discontinuities and satisfying symmetric jump (or
transmission) conditions. Their sampling result states that the
transformation%
\begin{equation*}
F\left( \lambda \right) =\int\limits_{0}^{\pi }f\left( x\right) u\left(
x,\lambda \right) dx,\text{ \ \ \ \ }f\in L^{2}\left( 0,\pi \right) ,
\end{equation*}%
can be reconstructed via the sampling form%
\begin{equation*}
F\left( \lambda \right) =\sum\limits_{n=0}^{\infty }F\left( \lambda
_{n}\right) \frac{\omega \left( \lambda \right) }{\left( \lambda -\lambda
_{n}\right) \omega ^{\prime }\left( \lambda _{n}\right) },
\end{equation*}%
where $\left\{ \lambda _{n}\right\} _{n=0}^{\infty }$ are the zeros of the
function $\omega \left( \lambda \right) ,$ and they are exactly the
eigenvalues of the discontinuous Sturm-Liouville problem. In this problem,
neither boundary conditions nor transmission conditions contain an
eigenparameter. In $\left[ 22\right] ,$ the author investigated the sampling
analysis associated with discontinuous Sturm-Liouville problems which has
one point of discontinuity and contains an eigenparameter in all boundary
conditions and derived sampling representations for transforms whose kernels
are either solutions or Green's functions.

We consider the boundary value problem:%
\begin{equation}
\tau \left( u\right) :=-u^{\shortparallel }\left( x\right) +q\left( x\right)
u\left( x\right) =\lambda u\left( x\right) ,~\ \ \ \ x\in I  \tag{1.9}
\end{equation}%
with boundary conditions;%
\begin{equation}
B_{1}\left( u\right) :=\beta _{1}u\left( a\right) +\beta _{2}u^{\prime
}\left( a\right) =0,  \tag{1.10}
\end{equation}%
\begin{equation}
B_{2}\left( u\right) :=\lambda \left( \alpha _{1}^{\prime }u\left( b\right)
-\alpha _{2}^{\prime }u^{\prime }\left( b\right) \right) -\left( \alpha
_{1}u\left( b\right) -\alpha _{2}u^{\prime }\left( b\right) \right) =0, 
\tag{1.11}
\end{equation}%
and transmission conditions at two points of discontinuity $c_{1}$and $%
c_{2}; $%
\begin{equation}
T_{1}\left( u\right) :=u\left( c_{1}^{-}\right) -\delta u\left(
c_{1}^{+}\right) =0,  \tag{1.12}
\end{equation}%
\begin{equation}
T_{2}\left( u\right) :=u^{\prime }\left( c_{1}^{-}\right) -\delta u^{\prime
}\left( c_{1}^{+}\right) +\lambda u\left( c_{1}^{-}\right) =0,  \tag{1.13}
\end{equation}%
\begin{equation}
T_{3}\left( u\right) :=\delta u\left( c_{2}^{-}\right) -\gamma u\left(
c_{2}^{+}\right) =0,  \tag{1.14}
\end{equation}%
\begin{equation}
T_{4}\left( u\right) :=\delta u^{\prime }\left( c_{2}^{-}\right) -\gamma
u^{\prime }\left( c_{2}^{+}\right) +\lambda u\left( c_{2}^{-}\right) =0, 
\tag{1.15}
\end{equation}%
where $I=[a,c_{1})\cup \left( c_{1},c_{2}\right) \cup
(c_{2},b],~a<c_{1}<c_{2}<b;$ $\lambda $ is a complex spectral parameter$;$ $%
q\left( x\right) $ is a given real valued function which is continuous in $%
[a,c_{1}),~\left( c_{1},c_{2}\right) $ and $(c_{2},b]$ and has finite limits 
$q(c_{1}^{\pm })=\underset{x\rightarrow c_{1}^{\pm }}{\lim }q\left( x\right) 
$ and $q(c_{2}^{\pm })=\underset{x\rightarrow c_{2}^{\pm }}{\lim }q\left(
x\right) ;$ $\alpha _{i},\alpha _{i}^{\prime },\beta _{i},$ $\delta ,\gamma
\in 
%TCIMACRO{\U{211d} }%
%BeginExpansion
\mathbb{R}
%EndExpansion
$ $\left( i=1,2\right) ,$ $\left\vert \beta _{1}\right\vert +\left\vert
\beta _{2}\right\vert \neq 0,$ $\delta >0,$ $\gamma \neq 0;$ $c_{i}^{\pm
}:=c_{i}\pm 0,$ $\left( i=1,2\right) $ and $\rho :=\left( \alpha
_{1}^{\prime }\alpha _{2}-\alpha _{1}\alpha _{2}^{\prime }\right) >0.$

In the present work, we investigate the sampling analysis associated with
discontinuous Sturm-Liouville problem (1.9)-(1.15) which has transmission
conditions at two points of discontinuity and contains an eigenparameter in
a boundary condition and two transmission conditions. This is the difference
between our problem and sampling theories associated with discontinuous
eigenproblems studied in the literature. The problem of deriving a sampling
theorem of Kramer type when the problem contains an eigenparameter in two
transmission conditions besides in a boundary conditions does not exist as
far as we know. To derive sampling theorem associated with the problem
(1.9)-(1.15), we study briefly the spectral properties of the problem
(1.9)-(1.15) and then we prove that integral transforms associated with the
problem (1.9)-(1.15) can also be reconstructed in a sampling form of
Lagrange interpolation type.

\section{Spectral Properties}

To formulate a theoretic approach to the problem (1.9)-(1.15) let $%
L:=L^{2}(a,c_{1})\oplus L^{2}(c_{1},c_{2})\oplus L^{2}\left( c_{2},b\right) $
and $%
%TCIMACRO{\U{2102} }%
%BeginExpansion
\mathbb{C}
%EndExpansion
^{3}:=%
%TCIMACRO{\U{2102} }%
%BeginExpansion
\mathbb{C}
%EndExpansion
\oplus 
%TCIMACRO{\U{2102} }%
%BeginExpansion
\mathbb{C}
%EndExpansion
\oplus 
%TCIMACRO{\U{2102} }%
%BeginExpansion
\mathbb{C}
%EndExpansion
$ and we define the Hilbert space $H=L\oplus 
%TCIMACRO{\U{2102} }%
%BeginExpansion
\mathbb{C}
%EndExpansion
^{3}$ with an inner product%
\begin{equation}
\begin{array}{l}
\left\langle F(.),G(.)\right\rangle _{H}:=\dint\limits_{a}^{c_{1}}f\left(
x\right) \overline{g}\left( x\right) dx+\delta
^{2}\dint\limits_{c_{1}}^{c_{2}}f\left( x\right) \overline{g}\left( x\right)
dx+\gamma ^{2}\dint\limits_{c_{2}}^{b}f\left( x\right) \overline{g}\left(
x\right) dx+ \\ 
\text{ \ \ \ \ \ \ \ \ \ \ \ \ \ \ \ \ \ \ \ \ \ \ }\dfrac{\gamma ^{2}}{\rho 
}h_{1}\overline{k_{1}}+h_{2}\overline{k_{2}}+\delta h_{3}\overline{k_{3}},%
\end{array}
\tag{2.1}
\end{equation}%
where

$F(x)=\left( 
\begin{array}{c}
f\left( x\right) \\ 
h_{1} \\ 
h_{2} \\ 
h_{3}%
\end{array}%
\right) ,G(x)=\left( 
\begin{array}{c}
g(x) \\ 
k_{1} \\ 
k_{2} \\ 
k_{3}%
\end{array}%
\right) \in H,$ $f\left( .\right) ,$ $g\left( .\right) \in L^{2}(a,b)$ and $%
h_{i},k_{i}\in 
%TCIMACRO{\U{2102} }%
%BeginExpansion
\mathbb{C}
%EndExpansion
$ $\left( i=1,2,3\right) .$

For convenience we put%
\begin{equation}
R_{b}\left( u\right) :=\alpha _{1}u\left( b\right) -\alpha _{2}u^{\prime
}\left( b\right) ,\text{ \ \ }R_{b}^{\prime }\left( u\right) :=\alpha
_{1}^{\prime }u\left( b\right) -\alpha _{2}^{\prime }u^{\prime }\left(
b\right) ,  \tag{2.2}
\end{equation}%
\begin{equation}
R_{c_{1}}\left( u\right) :=u\left( c_{1}^{-}\right) ,\text{ \ \ \ \ }%
R_{c_{1}}^{\prime }\left( u\right) :=u^{\prime }\left( c_{1}^{-}\right)
-\delta u^{\prime }\left( c_{1}^{+}\right) ,  \tag{2.3}
\end{equation}%
\begin{equation}
R_{c_{2}}\left( u\right) :=u\left( c_{2}^{-}\right) ,\text{ \ \ \ \ }%
R_{c_{2}}^{\prime }\left( u\right) :=\delta u^{\prime }\left(
c_{2}^{-}\right) -\gamma u^{\prime }\left( c_{2}^{+}\right) .  \tag{2.4}
\end{equation}

For functions $f(x),$ which is defined on $I$ and has finite limit $%
f(c_{1}^{\pm })=\underset{x\rightarrow c_{1}^{\pm }}{\lim }$ $f(x)$, $%
f(c_{2}^{\pm })=\underset{x\rightarrow c_{2}^{\pm }}{\lim }$ $f(x),$ by $%
f_{(i)}(x)$ $\left( i=\overline{1,4}\right) $ we denote the functions%
\begin{equation}
\begin{array}{l}
f_{(1)}(x):=\left\{ 
\begin{array}{l}
f(x),\text{ \ }x\in \lbrack a,c_{1}), \\ 
f(c_{1}^{-}),\text{ \ }x=c_{1},%
\end{array}%
\right. \ \ \ f_{(2)}(x):=\left\{ 
\begin{array}{l}
f(x),\text{ \ }x\in \left( c_{1},c_{2}\right) , \\ 
f(c_{1}^{+}),\text{ \ }x=c_{1},%
\end{array}%
\right. \\ 
\\ 
f_{(3)}(x):=\left\{ 
\begin{array}{l}
f(x),\text{ \ }x\in \left( c_{1},c_{2}\right) , \\ 
f(c_{2}^{-}),\text{ \ }x=c_{2},%
\end{array}%
\right. \text{ \ }f_{(4)}(x):=\left\{ 
\begin{array}{l}
f(x),\text{ \ }x\in \left( c_{2},b\right] , \\ 
f(c_{2}^{+}),\text{ \ }x=c_{2},%
\end{array}%
\right.%
\end{array}
\tag{2.5}
\end{equation}%
which are defined on $I_{1}:=[a,c_{1}],$ $I_{2}:=\left[ c_{1},c_{2}\right] $
and $I_{3}:=[c_{2},b],$ respectively.

In this space, we define a linear operator $A$ by the domain of definition%
\begin{equation}
D\left( A\right) :=\left\{ F\left( x\right) =\left( 
\begin{array}{c}
f\left( x\right) \\ 
R_{b}^{\prime }\left( f\right) \\ 
R_{c_{1}}\left( f\right) \\ 
R_{c_{2}}\left( f\right)%
\end{array}%
\right) \in H\left\vert 
\begin{array}{l}
f_{\left( i\right) }\left( .\right) ,\text{ }f_{\left( i\right) }^{\prime
}\left( .\right) \text{ are absolutely continuous } \\ 
\text{ in }I_{i}\text{ \ }\left( i=1,2,3\right) ;\text{ }\tau \left(
f\right) \in L\text{; } \\ 
B_{1}\left( f\right) =0,\text{ }T_{1}\left( f\right) =T_{3}\left( f\right)
=0, \\ 
h_{1}=R_{b}^{\prime }\left( f\right) ,h_{2}=R_{c_{1}}\left( f\right)
,h_{3}=R_{c_{2}}\left( f\right)%
\end{array}%
\right. \right\}  \tag{2.6}
\end{equation}%
and%
\begin{equation}
A\left( 
\begin{array}{c}
f\left( x\right) \\ 
R_{b}^{\prime }\left( f\right) \\ 
R_{c_{1}}\left( f\right) \\ 
R_{c_{2}}\left( f\right)%
\end{array}%
\right) =\left( 
\begin{array}{c}
\tau \left( f\right) \\ 
R_{b}\left( f\right) \\ 
-R_{c_{1}}^{\prime }\left( f\right) \\ 
-R_{c_{2}}^{\prime }\left( f\right)%
\end{array}%
\right) ,\text{ \ \ \ \ \ }\left( 
\begin{array}{c}
f\left( x\right) \\ 
R_{b}^{\prime }\left( f\right) \\ 
R_{c_{1}}\left( f\right) \\ 
R_{c_{2}}\left( f\right)%
\end{array}%
\right) \in D\left( A\right) .\text{\ \ }  \tag{2.7}
\end{equation}%
Consequently, the problem (1.9)-(1.15) can be rewritten in operator form as $%
AF=\lambda F,$ i.e., the problem (1.9)-(1.15) can be considered as the
eigenvalue problem for the operator $A.$

The following theorem can be proven by the same methods in similar studies $%
\left[ 21,23-25\right] .$

\begin{theorem}
\bigskip $\mathit{i)}$ $A$\textit{\ is a symmetric operator in }$H;$\textit{%
\ the eigenvalue problem for the operator }$A$\textit{\ and the problem
(1.9)-(1.15) coincide. }$\mathit{ii)}$\textit{\ All eigenvalues and
eigenfunctions of the operator }$A$ \textit{(or the problem (1.9)-(1.15))
are real. }$\mathit{iii)}$\textit{\ Two eigenfunctions }$u\left( x,\lambda
\right) $\textit{\ and }$v\left( x,\mu \right) $\textit{\ corresponding to
different eigenvalues }$\lambda $\textit{\ and }$\mu $\textit{, are
orthogonal, i.e.,}%
\begin{equation}
\begin{array}{l}
\dint\limits_{a}^{c_{1}}\mathit{\ }u\left( x,\lambda \right) v\left( x,\mu
\right) dx+\delta ^{2}\dint\limits_{c_{1}}^{c_{2}}u\left( x,\lambda \right)
v\left( x,\mu \right) dx+\gamma ^{2}\dint\limits_{c_{2}}^{b}u\left(
x,\lambda \right) v\left( x,\mu \right) dx+ \\ 
\text{\ \ \ \ \ }\dfrac{\gamma ^{2}}{\rho }R_{b}^{\prime }\left( u\right)
R_{b}^{\prime }\left( v\right) +R_{c_{1}}\left( u\right) R_{c_{1}}\left(
v\right) +\delta R_{c_{2}}\left( u\right) R_{c_{2}}\left( v\right) =0.%
\end{array}
\tag{2.8}
\end{equation}
\end{theorem}

Now, we will construct a special fundamental system of solutions of the
equation (1.9). By virtue of theorem 1.5 in $[26],$ we will define two
solutions of the equation (1.9) as follows:%
\begin{equation}
\phi _{\lambda }\left( x\right) =\left\{ 
\begin{array}{l}
\phi _{1\lambda }\left( x\right) ,\text{ \ \ }x\in \left[ a,c_{1}\right) ,
\\ 
\phi _{2\lambda }\left( x\right) ,\text{ \ }x\in \left( c_{1},c_{2}\right) ,
\\ 
\phi _{3\lambda }\left( x\right) ,\text{ \ \ }x\in \left( c_{2},b\right] ,%
\end{array}%
\right. ,\text{\ }\chi _{\lambda }\left( x\right) =\left\{ 
\begin{array}{l}
\chi _{1\lambda }\left( x\right) ,\text{ \ \ }x\in \left[ a,c_{1}\right) ,
\\ 
\chi _{2\lambda }\left( x\right) ,\text{ \ }x\in \left( c_{1},c_{2}\right) ,
\\ 
\chi _{3\lambda }\left( x\right) ,\text{ \ \ }x\in \left( c_{2},b\right] .%
\end{array}%
\right.  \tag{2.9}
\end{equation}

Let $\phi _{1\lambda }\left( x\right) =\phi _{1}\left( x,\lambda \right) $
be the solution of the equation (1.9) on $\left[ a,c_{1}\right] ,$ which
satisfies the initial conditions%
\begin{equation}
u(a)=\beta _{2},\text{ \ \ \ \ \ \ }u^{\prime }(a)=-\beta _{1}.  \tag{2.10}
\end{equation}

By virtue of theorem 1.5 in $[26],$ after defining this solution we may
define the solution $\phi _{2\lambda }\left( x\right) =\phi _{2}\left(
x,\lambda \right) $ of the equation (1.9) on $\left[ c_{1},c_{2}\right] $ by
means of the solution $\phi _{1\lambda }\left( x\right) $ by the nonstandard
initial conditions%
\begin{equation}
u(c_{1})=\frac{1}{\delta }\phi _{1\lambda }\left( c_{1}^{-}\right) ,\text{ \
\ \ \ \ }u^{\prime }(c_{1})=\frac{1}{\delta }\left\{ \phi _{1\lambda
}^{\prime }\left( c_{1}^{-}\right) +\lambda \phi _{1\lambda }\left(
c_{1}^{-}\right) \right\} .  \tag{2.11}
\end{equation}%
After defining this solution, we may define the solution $\phi _{3\lambda
}\left( x\right) =\phi _{3}\left( x,\lambda \right) $ of the equation (1.9)
on $\left[ c_{2},b\right] $ by means of the solution $\phi _{2\lambda
}\left( x\right) $ by the nonstandard initial conditions%
\begin{equation}
u(c_{2})=\frac{\delta }{\gamma }\phi _{2\lambda }\left( c_{2}^{-}\right) ,%
\text{ \ \ \ \ \ }u^{\prime }(c_{2})=\frac{1}{\gamma }\left\{ \delta \phi
_{2\lambda }^{\prime }\left( c_{2}^{-}\right) +\lambda \phi _{2\lambda
}\left( c_{2}^{-}\right) \right\} .  \tag{2.12}
\end{equation}%
Thus, $\phi _{\lambda }\left( x\right) =\phi \left( x,\lambda \right) $
satisfies the equation (1.9) on $I,$ yhe boundary condition (1.10) and the
transmission conditions (1.12)-(1.15).

Analogically, first we define the solution $\chi _{3\lambda }\left( x\right)
=\chi _{3}\left( x,\lambda \right) $ of the equation (1.9) on $\left[ c_{2},b%
\right] $ by the initial conditions%
\begin{equation}
u\left( b\right) =\lambda \alpha _{2}^{\prime }-\alpha _{2},\text{ \ \ \ \ \
\ \ }u^{\prime }\left( b\right) =\lambda \alpha _{1}^{\prime }-\alpha _{1}. 
\tag{2.13}
\end{equation}%
Again, after defining this solution, we define the solution $\chi _{2\lambda
}\left( x\right) =\chi _{2}\left( x,\lambda \right) $ of the equation (1.9)
on $\left[ c_{1},c_{2}\right] $ by the initial conditions%
\begin{equation}
u(c_{2})=\frac{\gamma }{\delta }\chi _{3\lambda }\left( c_{2}^{+}\right) ,%
\text{ \ \ \ \ \ }u^{\prime }(c_{2})=\frac{\gamma }{\delta }\left\{ \chi
_{3\lambda }^{\prime }\left( c_{2}^{+}\right) -\frac{\lambda }{\delta }\chi
_{3\lambda }\left( c_{2}^{+}\right) \right\} .  \tag{2.14}
\end{equation}%
After defining this solution, we define the solution $\chi _{1\lambda
}\left( x\right) =\chi _{1}\left( x,\lambda \right) $ of the equation (1.9)
on $\left[ a,c_{1}\right] $ by the initial conditions%
\begin{equation}
u(c_{1})=\delta \chi _{2\lambda }\left( c_{1}^{+}\right) ,\text{ \ \ \ \ \ }%
u^{\prime }(c_{1})=\delta \left\{ \chi _{2\lambda }^{\prime }\left(
c_{1}^{+}\right) -\lambda \chi _{2\lambda }\left( c_{1}^{+}\right) \right\} .
\tag{2.15}
\end{equation}%
Thus, $\chi _{\lambda }\left( x\right) =\chi \left( x,\lambda \right) $
satisfies the equation (1.9) on $I,$ the boundary condition (1.11) and the
transmission conditions (1.12)-(1.15).

Since the Wronskian $\ W\left( \phi _{i\lambda },\chi _{i\lambda };x\right) $
are independent on variable $x\in I_{i}$ $\ \left( i=1,2,3\right) $ and $%
\phi _{i\lambda }\left( x\right) $ and $\chi _{i\lambda }\left( x\right) $
are the entire functions of the parameter $\lambda $ for each $x\in I_{i},$
then the functions%
\begin{equation}
\omega _{i}\left( \lambda \right) :=W\left( \phi _{i\lambda },\chi
_{i\lambda };x\right) =\phi _{i\lambda }\left( x\right) \chi _{i\lambda
}^{\prime }\left( x\right) -\phi _{i\lambda }^{\prime }\left( x\right) \chi
_{i\lambda }\left( x\right) ,\text{ \ }\left( i=1,2,3\right)  \tag{2.16}
\end{equation}%
are the entire functions of parameter $\lambda .$ Taking into account
(2.11), (2.12), (2.14) and (2.15), after a short calculation, we get%
\begin{equation*}
\omega _{1}\left( \lambda \right) =\delta ^{2}\omega _{2}\left( \lambda
\right) =\gamma ^{2}\omega _{3}\left( \lambda \right)
\end{equation*}%
for each $\lambda \in 
%TCIMACRO{\U{2102} }%
%BeginExpansion
\mathbb{C}
%EndExpansion
.$

\begin{corollary}
The zeros of the functions $\omega _{1}\left( \lambda \right) ,$ $\omega
_{2}\left( \lambda \right) $ and $\omega _{3}\left( \lambda \right) $
coincide.
\end{corollary}

Then,\ we may introduce to the consideration the characteristic function $%
\omega \left( \lambda \right) $ as%
\begin{equation}
\omega \left( \lambda \right) :=\omega _{1}\left( \lambda \right) =\delta
^{2}\omega _{2}\left( \lambda \right) =\gamma ^{2}\omega _{3}\left( \lambda
\right) .  \tag{2.17}
\end{equation}

\begin{theorem}
The eigenvalues of the problem (1.9)-(1.15) are coincided zeros of the
function $\omega \left( \lambda \right) .$
\end{theorem}

\begin{proof}
...
\end{proof}

\begin{lemma}
All eigenvalues $\lambda _{n}$ of the problem (1.9)-(1.15) are simple zeros
of $\omega \left( \lambda \right) .$
\end{lemma}

\begin{proof}
...
\end{proof}

If $\lambda _{n}$ $\left( n=0,1,2,...\right) $ denote the zeros of $\omega
\left( \lambda \right) ,$ then 
\begin{equation}
\Phi _{\lambda _{n}}\left( x\right) :=\left( 
\begin{array}{c}
\phi _{\lambda _{n}}\left( x\right) \\ 
R_{b}^{\prime }\left( \phi _{\lambda _{n}}\right) \\ 
R_{c_{1}}\left( \phi _{\lambda _{n}}\right) \\ 
R_{c_{2}}\left( \phi _{\lambda _{n}}\right)%
\end{array}%
\right)  \tag{2.26}
\end{equation}%
are the corresponding eigenvectors of the operator $A,$ satisfying the
orthogonality relation%
\begin{equation}
\left\langle \Phi _{\lambda _{n}}\left( .\right) ,\Phi _{\lambda _{m}}\left(
.\right) \right\rangle _{H}=0\text{ \ for \ \ }n\neq m.  \tag{2.27}
\end{equation}

Here $\left\{ \phi _{\lambda _{n}}\left( .\right) \right\} _{n=0}^{\infty }$
is a sequence of eigenfunctions of the problem (1.9)-(1.15) corresponding to
the eigenvectors of $A,$ i.e.,%
\begin{equation}
\Psi _{\lambda _{n}}\left( x\right) :=\frac{\Phi _{\lambda _{n}}\left(
x\right) }{\left\Vert \Phi _{\lambda _{n}}\left( x\right) \right\Vert _{H}}%
=\left( 
\begin{array}{c}
\Psi _{\lambda _{n}}\left( x\right) \\ 
R_{b}^{\prime }\left( \Psi _{\lambda _{n}}\right) \\ 
R_{c_{1}}\left( \Psi _{\lambda _{n}}\right) \\ 
R_{c_{2}}\left( \Psi _{\lambda _{n}}\right)%
\end{array}%
\right) .  \tag{2.28}
\end{equation}

Therefore $\left\{ \chi _{\lambda _{n}}\left( .\right) \right\}
_{n=0}^{\infty }$ is another set of eigenfunctions which is related by $%
\left\{ \phi _{\lambda _{n}}\left( .\right) \right\} _{n=0}^{\infty }$ with%
\begin{equation}
\chi _{\lambda _{n}}\left( x\right) =k_{n}\phi _{\lambda _{n}}\left(
x\right) ,\text{ \ }x\in I,\text{ }n\in 
%TCIMACRO{\U{2124} }%
%BeginExpansion
\mathbb{Z}
%EndExpansion
\tag{2.29}
\end{equation}%
where $k_{n}\neq 0$ are non-zero constants, since all eigenvalues are simple.

\section{Asymptotic Formulas for Eigenvalues and Eigenfunctions}

The asymptotics formulas for eigenvalues and eigenfunctions can be derived
similar to the classical techniques of $\left[ 26\right] ,$ see also $\left[
23-25\right] .$ We state the results briefly.

Let $\phi _{\lambda }\left( .\right) $ be the solutions of the equation
(1.9) defined in section 2, and let $\lambda =s^{2}.$ Then the following
integral equations hold for $k=0$ and $k=1:$%
\begin{equation}
\begin{array}{l}
\dfrac{d^{k}}{dx^{k}}\phi _{1\lambda }\left( x\right) =\beta _{2}\dfrac{d^{k}%
}{dx^{k}}\left( \cos s\left( x-a\right) \right) -\dfrac{\beta _{1}}{s}\dfrac{%
d^{k}}{dx^{k}}(\sin s\left( x-a\right) )+ \\ 
\text{ \ \ \ \ \ \ \ \ \ \ \ \ \ \ \ \ \ \ \ \ }\dfrac{1}{s}\underset{a}{%
\overset{x}{\int }}\dfrac{d^{k}}{dx^{k}}\left( \sin s\left( x-y\right)
\right) q\left( y\right) \phi _{1\lambda }\left( y\right) dy,%
\end{array}
\tag{3.1}
\end{equation}%
\begin{equation}
\begin{array}{l}
\dfrac{d^{k}}{dx^{k}}\phi _{2\lambda }\left( x\right) =\dfrac{1}{\delta }%
\phi _{1\lambda }\left( c_{1}^{-}\right) \dfrac{d^{k}}{dx^{k}}\left( \cos
s\left( x-c_{1}\right) \right) +\dfrac{1}{s\delta }\left\{ \phi _{1\lambda
}^{\prime }\left( c_{1}^{-}\right) +s^{2}\phi _{1\lambda }\left(
c_{1}^{-}\right) \right\} \times \\ 
\text{ \ \ \ \ \ \ \ \ \ \ \ \ \ \ \ \ \ }\dfrac{d^{k}}{dx^{k}}\left( \sin
s\left( x-c_{1}\right) \right) +\dfrac{1}{s}\underset{c_{1}}{\overset{x}{%
\int }}\dfrac{d^{k}}{dx^{k}}\left( \sin s\left( x-y\right) \right) q\left(
y\right) \phi _{2\lambda }\left( y\right) dy,%
\end{array}
\tag{3.2}
\end{equation}%
\begin{equation}
\begin{array}{l}
\dfrac{d^{k}}{dx^{k}}\phi _{3\lambda }\left( x\right) =\dfrac{\delta }{%
\gamma }\phi _{2\lambda }\left( c_{2}^{-}\right) \dfrac{d^{k}}{dx^{k}}\left(
\cos s\left( x-c_{2}\right) \right) +\dfrac{1}{s\gamma }\left\{ \delta \phi
_{2\lambda }^{\prime }\left( c_{2}^{-}\right) +s^{2}\phi _{2\lambda }\left(
c_{2}^{-}\right) \right\} \times \\ 
\text{ \ \ \ \ \ \ \ \ \ \ \ \ \ \ \ \ \ \ }\dfrac{d^{k}}{dx^{k}}\left( \sin
s\left( x-c_{2}\right) \right) +\dfrac{1}{s}\underset{c_{2}}{\overset{x}{%
\int }}\dfrac{d^{k}}{dx^{k}}\left( \sin s\left( x-y\right) \right) q\left(
y\right) \phi _{3\lambda }\left( y\right) dy.%
\end{array}
\tag{3.3}
\end{equation}

Let $\lambda =s^{2}$ and $\left\vert \func{Im}s\right\vert =t.$ Then the
functions $\phi _{i\lambda }\left( x\right) $ have the following asymptotic
representations for $\left\vert \lambda \right\vert \rightarrow \infty ,$
which hold uniformly for $x\in I_{i}$ $\left( i=1,2,3\right) :$%
\begin{equation}
\dfrac{d^{k}}{dx^{k}}\phi _{1\lambda }\left( x\right) =\beta _{2}\dfrac{d^{k}%
}{dx^{k}}\left( \cos s\left( x-a\right) \right) +O\left( \frac{1}{s}%
e^{t\left( x-a\right) }\right) ,  \tag{3.4}
\end{equation}%
\begin{equation}
\dfrac{d^{k}}{dx^{k}}\phi _{2\lambda }\left( x\right) =\dfrac{s\beta _{2}}{%
\delta }\cos \left( s\left( c_{1}-a\right) \right) \dfrac{d^{k}}{dx^{k}}%
\left( \sin s\left( x-c_{1}\right) \right) +O\left( e^{t\left( x-a\right)
}\right) ,  \tag{3.5}
\end{equation}%
\begin{equation}
\dfrac{d^{k}}{dx^{k}}\phi _{3\lambda }\left( x\right) =\dfrac{s^{2}\beta _{2}%
}{\delta \gamma }\cos \left( s\left( c_{1}-a\right) \right) \sin \left(
s\left( c_{2}-c_{1}\right) \right) \dfrac{d^{k}}{dx^{k}}\left( \sin s\left(
x-c_{2}\right) \right) +O\left( se^{t\left( x-a\right) }\right) ,  \tag{3.6}
\end{equation}

By substituting (3.6) into the representation%
\begin{equation}
\omega \left( \lambda \right) =\gamma ^{2}\left\{ \left( s^{2}\alpha
_{1}^{\prime }-\alpha _{1}\right) \phi _{3\lambda }\left( b\right) -\left(
s^{2}\alpha _{2}^{\prime }-\alpha _{2}\right) \phi _{3\lambda }^{\prime
}\left( b\right) \right\} ,  \tag{3.7}
\end{equation}%
then the characteristic function $\omega \left( \lambda \right) $ has the
following asymptotic representation:%
\begin{equation}
\omega \left( \lambda \right) =\frac{s^{5}\gamma \beta _{2}\alpha
_{2}^{\prime }}{\delta }\cos \left( s\left( c_{1}-a\right) \right) \sin
\left( s\left( c_{2}-c_{1}\right) \right) \cos \left( s\left( b-c_{2}\right)
\right) +O\left( s^{4}e^{t\left( b-a\right) }\right) .  \tag{3.8}
\end{equation}

\begin{corollary}
The eigenvalues of the problem (1.9)-(1.15) are bounded below.
\end{corollary}

Now we can obtain the asymptotic approximation formula for the eigenvalues
of the problem (1.9)-(1.15). Since the eigenvalues coincide with the zeros
of the entire function $\omega \left( \lambda \right) ,$ it follows that
they have no finite limit. Moreover, we know that all eigenvalues are real
and bounded below. Therefore, we may renumber them as $\lambda _{0}\leq
\lambda _{1}\leq ...,$ listed according to their multiplicity.

\begin{theorem}
The eigenvalues $\lambda _{n}=s_{n}^{2},\left( n=0,1,...\right) $ of the
problem (1.9)-(1.15) have the following asymptotic representation for $%
n\rightarrow \infty :$%
\begin{equation}
\widetilde{s}_{n}=\frac{\left( n+1/2\right) \pi }{\left( c_{1}-a\right) }%
+O\left( \frac{1}{n}\right) ,  \tag{3.9}
\end{equation}%
\begin{equation}
\widetilde{\widetilde{s}}_{n}=\frac{n\pi }{\left( c_{2}-c_{1}\right) }%
+O\left( \frac{1}{n}\right) ,  \tag{3.10}
\end{equation}%
\begin{equation}
\widetilde{\widetilde{\widetilde{s}}}_{n}=\frac{\left( n+1/2\right) \pi }{%
\left( b-c_{2}\right) }+O\left( \frac{1}{n}\right) .  \tag{3.11}
\end{equation}
\end{theorem}

\begin{proof}
...
\end{proof}

Then from (3.4)-(3.6) (for $k=0$) and the above theorem, the asymptotic
behaviour of the eigenfunctions:%
\begin{equation}
\phi _{\lambda _{n}}\left( x\right) =\left\{ 
\begin{array}{l}
\phi _{1\lambda _{n}}\left( x\right) ,\text{ \ \ }x\in \left[ a,c_{1}\right)
, \\ 
\phi _{2\lambda _{n}}\left( x\right) ,\text{ \ }x\in \left(
c_{1},c_{2}\right) , \\ 
\phi _{3\lambda _{n}}\left( x\right) ,\text{ \ \ }x\in \left( c_{2},b\right]
,%
\end{array}%
\right.  \tag{3.13}
\end{equation}%
of the problem (1.9)-(1.15) is given by%
\begin{equation*}
\phi _{\widetilde{\lambda }_{n}}\left( x\right) =\left\{ 
\begin{array}{l}
\beta _{2}\cos \left( \dfrac{\left( n+1/2\right) \pi }{\left( c_{1}-a\right) 
}\left( x-a\right) \right) +O\left( \dfrac{1}{n}\right) ,\text{ \ \ }x\in %
\left[ a,c_{1}\right) , \\ 
O\left( 1\right) ,\text{ \ \ \ \ \ \ \ \ \ \ \ \ \ \ \ \ \ \ \ \ \ \ \ \ \ \
\ \ \ \ \ \ \ \ \ \ \ \ \ \ \ \ \ \ \ }x\in \left( c_{1},c_{2}\right) , \\ 
O\left( 1\right) ,\text{ \ \ \ \ \ \ \ \ \ \ \ \ \ \ \ \ \ \ \ \ \ \ \ \ \ \
\ \ \ \ \ \ \ \ \ \ \ \ \ \ \ \ \ \ \ }x\in \left( c_{2},b\right] ,%
\end{array}%
\right.
\end{equation*}%
\begin{equation*}
\phi _{\widetilde{\widetilde{\lambda }}_{n}}\left( x\right) =\left\{ 
\begin{array}{l}
\beta _{2}\cos \left( \dfrac{n\pi }{\left( c_{2}-c_{1}\right) }\left(
x-a\right) \right) +O\left( \dfrac{1}{n}\right) ,\text{ \ \ \ \ \ \ \ \ \ \
\ \ \ \ \ \ \ \ \ \ \ \ \ \ \ }x\in \left[ a,c_{1}\right) , \\ 
\dfrac{n\pi }{\left( c_{2}-c_{1}\right) }\dfrac{\beta _{2}}{\delta }\cos
\left( \dfrac{n\pi \left( c_{1}-a\right) }{\left( c_{2}-c_{1}\right) }%
\right) \sin \left( \dfrac{n\pi }{\left( c_{2}-c_{1}\right) }\left(
x-c_{1}\right) \right) +O\left( 1\right) ,\text{\ }x\in \left(
c_{1},c_{2}\right) , \\ 
O\left( 1\right) ,\text{ \ \ \ \ \ \ \ \ \ \ \ \ \ \ \ \ \ \ \ \ \ \ \ \ \ \
\ \ \ \ \ \ \ \ \ \ \ \ \ \ \ \ \ \ \ \ \ \ \ \ \ \ \ \ \ \ \ \ \ \ \ \ \ \
\ \ \ }x\in \left( c_{2},b\right] ,%
\end{array}%
\right.
\end{equation*}%
\begin{equation*}
\phi _{\widetilde{\widetilde{\widetilde{\lambda }}}_{n}}\left( x\right)
=\left\{ 
\begin{array}{l}
\beta _{2}\cos \left( \dfrac{\left( n+1/2\right) \pi }{\left( b-c_{2}\right) 
}\left( x-a\right) \right) +O\left( \dfrac{1}{n}\right) ,\text{ \ \ \ \ \ \
\ \ \ \ \ \ \ \ \ \ \ \ \ \ \ \ \ \ }x\in \left[ a,c_{1}\right) , \\ 
\dfrac{\left( n+1/2\right) \pi }{\left( b-c_{2}\right) }\dfrac{\beta _{2}}{%
\delta }\cos \left( \dfrac{\left( n+1/2\right) \pi \left( c_{1}-a\right) }{%
\left( b-c_{2}\right) }\right) \sin \left( \dfrac{\left( n+1/2\right) \pi
\left( c_{2}-c_{1}\right) }{\left( b-c_{2}\right) }\right) \times \\ 
\sin \left( \dfrac{\left( n+1/2\right) \pi }{\left( b-c_{2}\right) }\left(
x-c_{2}\right) \right) +O\left( 1\right) ,\text{\ \ \ \ \ \ \ \ \ \ \ \ \ \
\ \ \ \ \ \ \ \ \ \ \ \ \ \ \ \ \ \ \ }x\in \left( c_{1},c_{2}\right) , \\ 
O\left( 1\right) ,\text{ \ \ \ \ \ \ \ \ \ \ \ \ \ \ \ \ \ \ \ \ \ \ \ \ \ \
\ \ \ \ \ \ \ \ \ \ \ \ \ \ \ \ \ \ \ \ \ \ \ \ \ \ \ \ \ \ \ \ \ \ \ \ \ \
\ \ \ \ \ \ \ \ }x\in \left( c_{2},b\right] ,%
\end{array}%
\right.
\end{equation*}%
all these asymptotic formulas hold uniformly for $x.$

\section{The Sampling Theorem}

Now we can derive the sampling theorem associated with the problem
(1.9)-(1.15).

\begin{theorem}
Let%
\begin{equation}
\phi _{\lambda }\left( x\right) =\left\{ 
\begin{array}{l}
\phi _{1\lambda }\left( x\right) ,\text{ \ \ }x\in \left[ a,c_{1}\right) ,
\\ 
\phi _{2\lambda }\left( x\right) ,\text{ \ }x\in \left( c_{1},c_{2}\right) ,
\\ 
\phi _{3\lambda }\left( x\right) ,\text{ \ \ }x\in \left( c_{2},b\right] ,%
\end{array}%
\right.  \tag{4.1}
\end{equation}%
be the solution of the equation (1.9) together with the conditions
(2.10)-(2.12). Let $g\left( .\right) \in L^{2}\left( a,b\right) $ and 
\begin{equation}
\digamma \left( \lambda \right) =\dint\limits_{a}^{c_{1}}\mathit{\ g}\left(
x\right) \phi _{1\lambda }\left( x\right) dx+\delta
^{2}\dint\limits_{c_{1}}^{c_{2}}\mathit{g}\left( x\right) \phi _{2\lambda
}\left( x\right) dx+\gamma ^{2}\dint\limits_{c_{2}}^{b}\mathit{g}\left(
x\right) \phi _{3\lambda }\left( x\right) dx.  \tag{4.2}
\end{equation}%
Then $\digamma \left( \lambda \right) $ is an entire function of exponential
type that can be reconstructed from its values at the points $\left\{
\lambda _{n}\right\} _{n=0}^{\infty }$ via the sampling formula 
\begin{equation}
\digamma \left( \lambda \right) =\sum_{n=0}^{\infty }\digamma \left( \lambda
_{n}\right) \frac{\omega \left( \lambda \right) }{\left( \lambda -\lambda
_{n}\right) \omega ^{\prime }\left( \lambda _{n}\right) }.  \tag{4.3}
\end{equation}%
The series (4.3) converges absolutely on $%
%TCIMACRO{\U{2102} }%
%BeginExpansion
\mathbb{C}
%EndExpansion
$ and uniformly on compact subset of $%
%TCIMACRO{\U{2102} }%
%BeginExpansion
\mathbb{C}
%EndExpansion
.$ Here $\omega \left( \lambda \right) $ is the entire function defined in
(2.17).
\end{theorem}

\begin{proof}
...
\end{proof}

\begin{remark}
To see that expansion (4.3) is a Lagrange type interpolation, we may replace 
$\omega \left( \lambda \right) $ by the canonical product%
\begin{equation}
\varpi \left( \lambda \right) =\left\{ 
\begin{array}{c}
\overset{\infty }{\underset{n=0}{\Pi }}\left( 1-\dfrac{\lambda }{\lambda _{n}%
}\right) ,\text{ \ if zero is not an eigenvalue,} \\ 
\\ 
\lambda \overset{\infty }{\underset{n=1}{\Pi }}\left( 1-\dfrac{\lambda }{%
\lambda _{n}}\right) ,\text{ \ if }\lambda _{0}=0\text{ is an eigenvalue.}%
\end{array}%
\right.  \tag{4.26}
\end{equation}
\end{remark}

From Hadamard's factorization theorem, see $\left[ 6\right] ,\omega \left(
\lambda \right) =h\left( \lambda \right) \varpi \left( \lambda \right) ,$
where $h\left( \lambda \right) $ is an entire function with no zeros. Thus,%
\begin{equation}
\frac{\omega \left( \lambda \right) }{\omega ^{\prime }\left( \lambda
_{n}\right) }=\frac{h\left( \lambda \right) \varpi \left( \lambda \right) }{%
h\left( \lambda _{n}\right) \varpi ^{\prime }\left( \lambda _{n}\right) } 
\tag{4.27}
\end{equation}%
and (4.2), (4.3) remain valid for the function $\digamma \left( \lambda
\right) /h\left( \lambda \right) .$ Hence%
\begin{equation}
\digamma \left( \lambda \right) =\overset{\infty }{\underset{n=0}{\sum }}%
\digamma \left( \lambda _{n}\right) \frac{h\left( \lambda \right) \varpi
\left( \lambda \right) }{\left( \lambda -\lambda _{n}\right) h\left( \lambda
_{n}\right) \varpi ^{\prime }\left( \lambda _{n}\right) }.  \tag{4.28}
\end{equation}%
We may redefine (4.2) by taking kernel $\phi _{\lambda }\left( .\right)
/h\left( \lambda \right) =\overline{\phi }_{\lambda }\left( .\right) $ to get%
\begin{equation}
\overline{\digamma }\left( \lambda \right) =\frac{\digamma \left( \lambda
\right) }{h\left( \lambda \right) }=\overset{\infty }{\underset{n=0}{\sum }}%
\overline{\digamma }\left( \lambda _{n}\right) \frac{\varpi \left( \lambda
\right) }{\left( \lambda -\lambda _{n}\right) \varpi ^{\prime }\left(
\lambda _{n}\right) }.  \tag{4.29}
\end{equation}

\end{document}